\documentclass[12pt]{article}
\usepackage{mathrsfs}
\usepackage{amsmath}
\usepackage{amsthm}
\usepackage{bbm}
\usepackage{subfloat}
\usepackage{subfigure}
\usepackage{fancyhdr,graphicx}
\usepackage{amsfonts}
\usepackage{amssymb}
\usepackage{latexsym,bm}
\usepackage{newlfont}
\usepackage{float}
\allowdisplaybreaks[3]
\usepackage[numbers,sort&compress]{natbib}
\usepackage{caption}

\title{{Counterexamples to a conjecture of Balasubramanian and Parthasarathy\thanks{This
work is supported by NSFC (grant no. 11401044).}}}
\author{Shunyi Liu\footnote{E-mail address: liu@chd.edu.cn.}
\\\small{College of Science, Chang'an University, Xi'an, Shaanxi 710064, P.R. China}}

\date{}

\begin{document}

\maketitle

\makeatletter
\newcommand{\rmnum}[1]{\romannumeral #1} ¡¡¡¡
\newcommand{\Rmnum}[1]{\expandafter\@slowromancap\romannumeral #1@}
\makeatother

\begin{abstract}
In 1980, Balasubramanian and Parthasarathy introduced the bivariate permanent polynomials of graphs and conjectured that this graph polynomial is a graph characterising polynomial, that is, any two graphs with the same bivariate permanent polynomial are isomorphic. In this note, we give counter examples to this conjecture.


\end{abstract}

\section{Introduction}
\setlength{\unitlength}{1cm}
A graph invariant is a function $f$ from the set of all graphs $\mathscr{G}$ into any commutative ring $R$ such that $f$ takes the same value on isomorphic graphs. When $R$ is a ring of polynomials in one or more variables, the invariant $f$ is called an invariant polynomial for graphs (or a \emph{graph polynomial}). As a graph invariant $f$ can be used to check whether two graphs are not isomorphic. If a graph polynomial $f$ also satisfies the converse condition that $f(G)=f(H)$ implies $G$ and $H$ are isomorphic, then $f$ is called a \emph{graph characterising polynomial}.

Many graph polynomials have been defined and extensively studied, like the characteristic, chromatic, matching, and Tutte polynomials (see, for example, \cite{ElMe1,ElMe2,Gut}). Besides their intrinsic interest, they encode useful combinatorial information about the given graph. In general, graph polynomials have been developed for measuring combinatorial graph invariants and for characterizing graphs. The latter is related to graph isomorphism problem and it is of interest to determine its ability to characterize graphs for any graph polynomial \cite{Noy}. One might ask whether or not we can find a graph characterising polynomial. To date, no useful graph characterising polynomials have been found. Indeed, all the graph polynomials mentioned above are not graph characterising polynomials.

In 1980, Balasubramanian and Parthasarathy~\cite{BaPa,Par} introduced a graph polynomial, bivariate permanent polynomial, and conjectured that this polynomial is a graph characterising polynomial. As far as I know, this conjecture is still open. In what follows, we shall call this conjecture the bivariate permanent polynomial conjecture (BPPC for short).

The \emph{permanent} of an $n\times n$ matrix $M$ with entries
$m_{ij}$ $(i,j=1, 2, \dots, n)$ is defined by
\begin{equation*}
\mathrm{per}(M)=\sum_{\sigma}\prod_{i=1}^{n}m_{i\sigma(i)},
\end{equation*}
where the sum is taken over all permutations $\sigma$ of $\{1, 2, \dots, n\}$. This scalar function of the matrix $M$ appears repeatedly in the literature of combinatorics and graph theory in connection with certain enumeration and extremal problems. For example, the permanent of a (0,1)-matrix enumerates perfect matchings in bipartite graphs \cite{LoPl}. The permanent is defined similarly to the determinant. However, no efficient algorithm for computing the permanent is known, while the determinant can be calculated using Gaussian elimination. More precisely, Valiant~\cite{Val} has shown that computing the permanent is $\#$P-complete even when restricted to (0,1)-matrices.

Let $G$ be a graph on $n$ vertices and $\bar{G}$ the complement of $G$. We use $A$ and $\bar{A}$ to denote the adjacency matrices
of $G$ and $\bar{G}$, respectively. The \emph{bivariate permanent polynomial}~\cite{BaPa} of $G$, $P(G;x,\lambda)$, is defined as
\begin{equation*}
P(G;x,\lambda)=\mathrm{per}(xI_n+\lambda A+\bar{A}),
\end{equation*}
where $I_n$ is the identity matrix of size $n$. Two graphs $G$ and $H$ are called \emph{copermanent} if they have the same bivariate permanent polynomial. A graph $H$, copermanent but non-isomorphic to a graph $G$, is called a \emph{copermanent mate} of $G$. We say that a graph $G$ is \emph{characterized} by its bivariate permanent polynomial if it has no copermanent mates. Thus BPPC can be restated that each graph is characterized by its bivariate permanent polynomial.

The validity of BPPC has been verified for all graphs on at most 7 vertices~\cite{BaPa}. In \cite{Par}, Parthasarathy showed that BPPC implies the celebrated graph reconstruction conjectures and gave a possible approach to settle BPPC. We were surprised when our literature search turned up only these two articles on the bivariate permanent polynomial. This may be due to the difficulty of actually computing the permanent.

It is worth pointing out that a univariate graph polynomial related to the matrix function permanent, named \emph{permanental polynomial}, has also been introduced by Merris et al.~\cite{MeRW}. The permanental polynomial of a graph $G$ is defined as $\mathrm{per}(xI_n-A)$, where $A$ is the adjacency matrix of $G$. It should be noted that the permanental polynomial is not a graph characterising polynomial~\cite{MeRW, Bor}. Characterizing graphs using permanental polynomial has recently been studied (see, for exmaple, \cite{LiZh1, LiZh2, ZhWL}).

In this paper, we give counter examples to BPPC by a computer search technique. More specifically, we determine the bivariate permanent polynomials for all graphs on at most 10 vertices, and count the number of graphs for which there exists at least one copermanent mate.

\section{Counterexamples}
To determine the bivariate permanent polynomials of graphs we first of all have to generate the graphs by computer. All graphs on at most 10 vertices are generated by the well-known nauty and Traces package~\cite{McPi}. Then the bivariate permanent polynomials of these graphs are computed by a Maple procedure. Finally we count the number of copermanent graphs.

The results are in Table~\ref{tab:Table 1}. This table lists for $n\le 10$ the total number of graphs on $n$ vertices, the total number of distinct bivariate permanent polynomials of such graphs, the number of such graphs with a copermanent mate, the fraction of such graphs with a copermanent mate, and the size of the largest family of copermanent graphs.

\begin{table}[htb]
\captionsetup{singlelinecheck=off,skip=0pt}
\caption{Computational data on $n\le 10$ vertices}
\label{tab:Table 1}
    \begin{tabular}{rrrrrr}\hline
    $n$ & $\#$graphs & $\#$perm. pols & $\#$ with coperm. mate & frac. with mate & max. family\\ \hline
     0  & 1          & 1              & 0                      &0                & 1 \\
     1  & 1          & 1              & 0                      &0                & 1 \\
     2  & 2          & 2              & 0                      &0                & 1 \\
     3  & 4          & 4              & 0                      &0                & 1 \\
     4  & 11         & 11             & 0                      &0                & 1 \\
     5  & 34         & 34             & 0                      &0                & 1 \\
     6  & 156        & 156            & 0                      &0                & 1 \\
     7  & 1044       & 1044           & 0                      &0                & 1 \\
     8  & 12346      & 12344          & 4                      & 0.000324         & 2 \\
     9  & 274668     & 274624         & 88                     & 0.000320         & 2 \\
     10 & 12005168   & 12004460       & 1416                   & 0.000118         & 2 \\ \hline
    \end{tabular}
\end{table}

In Table~\ref{tab:Table 1} we see that there are 4 graphs on eight vertices, 88 graphs on nine vertices, and 1416 graphs on ten vertices are not characterized by their bivariate permanent polynomials. Although the bivariate permanent polynomial is not a graph characterising polynomial, Table~\ref{tab:Table 1} gives some indication that possibly the fraction of graphs with a copermanent mate tends to zero as $n$ tends to infinity.

\begin{figure}[tphb]
    \begin{center}
\includegraphics[scale=0.65]{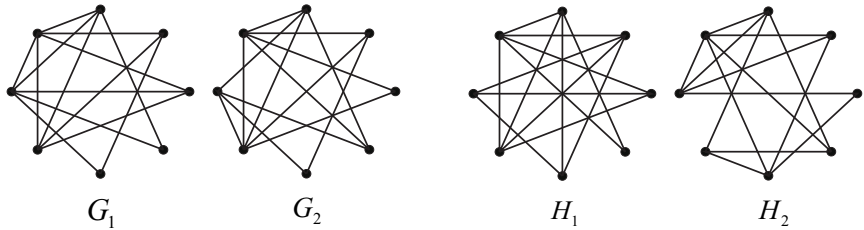}
\end{center}
   \caption{Two pairs of copermanent graphs on 8 vertices.}
   \protect\label{fig:Figure_1}
\end{figure}

Two pairs of copermanent graphs on 8 vertices are given in Figure~\ref{fig:Figure_1}, where $G_1$ and $G_2$ are copermanent, and $H_1$ and $H_2$ are copermanent. Their bivariate permanent polynomials are
\begin{align*}
&\quad P(G_1;x,\lambda)=P(G_2;x,\lambda)\\
&=x^8+14x^6\lambda^2+14x^6+12x^5\lambda^3+44x^5\lambda^2+44x^5\lambda+12x^5+69x^4\lambda^4+112x^4\lambda^3+268x^4\lambda^2\\
&\quad+112x^4\lambda+69x^4+82x^3\lambda^5+402x^3\lambda^4+748x^3\lambda^3+748x^3\lambda^2+402x^3\lambda+82x^3+130x^2\lambda^6\\
&\quad+648x^2\lambda^5+1804x^2\lambda^4+2256x^2\lambda^3+1804x^2\lambda^2+648x^2\lambda+130x^2+88x\lambda^7+742x\lambda^6\\
&\quad+2434x\lambda^5+4152x\lambda^4+4152x\lambda^3+2434x\lambda^2+742x\lambda+88x+40\lambda^8+360\lambda^7+1520\lambda^6\\
&\quad+3320\lambda^5+4353\lambda^4+3320\lambda^3+1520\lambda^2+360\lambda+40,
\end{align*}
and
\begin{align*}
&\quad P(H_1;x,\lambda)=P(H_2;x,\lambda)\\
&=x^8+14x^6\lambda^2+14x^6+10x^5\lambda^3+46x^5\lambda^2+46x^5\lambda+10x^5+69x^4\lambda^4+108x^4\lambda^3+276x^4\lambda^2\\
&\quad+108x^4\lambda+69x^4+78x^3\lambda^5+418x^3\lambda^4+736x^3\lambda^3+736x^3\lambda^2+418x^3\lambda+78x^3+144x^2\lambda^6\\
&\quad+672x^2\lambda^5+1814x^2\lambda^4+2160x^2\lambda^3+1814x^2\lambda^2+672x^2\lambda+144x^2+130x\lambda^7+830x\lambda^6\\
&\quad+2412x\lambda^5+4044x\lambda^4+4044x\lambda^3+2412x\lambda^2+830x\lambda+130x+52\lambda^8+468\lambda^7+1672\lambda^6\\
&\quad+3208\lambda^5+4033\lambda^4+3208\lambda^3+1672\lambda^2+468\lambda+52.
\end{align*}


\end{document}